# POLYNOMIAL RINGS WHOSE PRIMES ARE SET THEORETIC COMPLETE INTERSECTIONS

## VAHAP ERDOĞDU

*ABSTRACT. A well known cojecture states that over a field K, height two prime ideals of the polynomial ring K[X, Y, Z] are set theoretic complete intersections. Motivated by this here in this note we seek for conditions on R over which the polynomial ring R[X] has the property that all its prime ideals are set theoretic complete intersections.*

## 0. Introduction

The question of relating the number of generators of an ideal to the height of the ideal has been, and still is, one of the interesting research topics in commutative algebra. This question, apparently, was first considered by Kronecker in late $19^{th}$ century. Since then an enormous amount of research has evolved around these types of questions (see, e.g. [1], [2], [5, Chapter V], [6] and the references there in). Among them still remaining unsolved is the question whether height two ideals in K[X, Y, Z] are set theoretic complete intersections, where K is of characteristic zero. When K is of positive characteristic it is known that height n-1 ideals in the polynomial rings in n variables over K are set theoretic complete intersections [1]. The information provided in this note may be of some help towards the solution of the conjecture when K is of characteristic zero.

## 1. Statements and their motivations

Let R be a commutative ring with identity and R[X] be the polynomial ring over R. One of our goals here is to prove the following:

**Theorem 1.** *Let R be a Noetherian integral domain of Krull dimension one which contains a field of characteristic zero. Then each prime ideal of R[X] is a set theoretic complete intersection if and only if R is a Dedekind domain with torsion ideal class group.*

The above statement is proved in Theorem 2.1 of [3] under the assumption that R is normal. The reason for replacing the normality by the condition that R contains a field of *characteristic* zero is that in Remark 2.2 of [3] it was pointed out that if K is a field of characteristic zero and R is the sub ring of K[[Y]] consisting of those power series whose Y term has coefficient zero, then R is not normal and that there is a height one prime ideal of R[X] that is not a set theoretic complete intersection. However, if K is of characteristic $p$ ($p$ a positive prime), then an argument similar to the following one will show that each prime ideal of R[X] is a set theoretic complete intersection. Let K be a field of characteristic $p$ and R = K $[Y^2, Y^3]$, then R is a Noetherian integral domain of Krull dimension one whose normalization S = K[Y] is a principal ideal domain, and therefore by Theorem 2.1 of [3], each prime ideal of S[X] is a set theoretic complete intersection. Let now P be any height one prime ideal of R[X]. Then $rad$(PS[X]) is the intersection of finite number of height one primes in S[X]. Since S[X] is a U.F.D., height one primes in S[X] are principal and that the intersections of a finite number of principal ideals in S[X] is again a principal ideal, and so it follows that $rad$(PS[X]) = $rad$(f) for some polynomial f in S[X]. Now K being of characteristic $p$, $f^p$ is in R[X] and therefore is in P. (This is because the *p-th* power of every element of S[X] is in R[X]). But then it is clear that P = $rad(f^p)$ in R[X]. Thus each height one prime in R[X] (in fact each prime in R[X]) is a set theoretic complete intersection without R being normal.

These facts give raise to the following result:

**Theorem 2.** *Let R be an integral domain (not necessarily Noetherian) containing a field of characteristic zero. If each height one prime ideal of R[X] is a set theoretic complete intersection, then R is normal.*

## 2. Proofs of the statements

Before proving the statements we recall that an ideal I of a ring R is a set theoretic complete intersection if there exist n elements $a_1, a_2, \ldots, a_n$ in R such that $rad($ I $) = rad(a_1, a_2, \ldots, a_n)$ where n = height of I. We first prove Theorem 2 and use it to prove Theorem 1.

*Proof of Theorem 2.* Let $R'$ be the integral closure of R in K, the field of fractions of R, and let $a$ be any element of $R'$. Then (X − a)K[X] is a prime ideal of K[X] and so (X − a)K[X] ∩ $R'$[X] = (X − a) $R'$[X] is a height one prime ideal of $R'$[X]. Therefore P = (X − a) $R'$[X] ∩ R[X] is a height one prime ideal of R[X] and hence by the hypothesis P = $rad($ f ), for some polynomial f in P. Since (X − a) $R'$[X] ∩ $R' = 0$, there is a prime Q in K[X] such that Q ∩ $R'$[X] = (X − a) $R'$[X] and so it follows that Q∩R[X] = P. Thus (X − a)K[X] and Q are two primes in K[X] lying over P. Since P∩R = 0, there can be only one prime in K[X] lying over P. Therefore PK[X] = $rad($ f )K[X] = (X − a)K[X]. It is now clear that in the U.F.D. K[X], the only prime divisor of f is (X − a). Hence it follows that f = $c(X - a)^n$ for some positive integer n and some constant $c$ in K. We now show that $c$ has to be the identity of R. For this we consider the ring homomorphism $\phi$ : R[X] → R[a] given by $\phi($ X $) = a$. As R[a] is an integral domain, we see that Ker$\phi$ is a prime ideal of R[X]. Moreover we have Ker$\phi$ = P in R[X]. This is because if g is any element of P, then g is in $PR'[X] \subseteq (X-a)R'[X]$ and so $g(a) = 0$ which implies that $P \subseteq \text{Ker}\,\phi$. Now using the fact that $\text{Ker}\,\phi \cap R = 0$ and the fact that there can not be a chain of three distinct primes in R[X] contracting to 0 in R, we obtain $P = \text{Ker}\,\phi$. Since $P = rad(f) = \text{Ker}\,\phi$ contains the monic polynomial for $a$, we see that f is monic and therefore $c = 1$ in R. Thus $f = (X-a)^n$ and so by comparing the coefficients we see that $na \in R$. Now using the fact that n is invertible in R we get $a \in R$, showing that R is normal.

Now that Theorem 2 is proved we prove Theorem1.

*Proof of Theorem 1.* If each prime in R[X] is a set theoretic complete intersection, then it follows from Theorem 2 that R is a Noetherian normal domain of dimension one and hence is a Dedekind domain. Let now P be any non zero prime ideal of R. Then as R is Noetherian, PR[X] is a height one prime ideal of R[X] and so by the assumption $PR[X] = rad(f)$ for some polynomial f in R[X]. Let $a$ be any element of P. Then $a$ is in $PR[X] = rad(f)$ and hence there is a positive integer n such that $a^n \in (f)$. Thus $a^n = fg$ for some polynomial g in R[X]. Since R[X] is an integral domain, by comparing the degrees of the polynomials on both sides, we see that $f = b$ for some constant $b$ in P. Thus $P = rad(b)$ in R and hence $P^m = Rb$ for some positive integer m. But then it follows that the class group of R is torsion.
Conversely, if R is a Dedekind domain with torsion class group, then by Theorem 2.1 of [1] each prime ideal of R[X] is a set theoretic complete intersection.

**Remark 3.** Let K be a field of characteristic zero and $R = K[X, Y]$. Then R is a polynomial ring in Y over the principal ideal domain K[X] and hence it follows from Theorem 1 that each prime ideal of R is a set theoretic complete intersection. Let P* be any height two prime ideal of $K[X, Y, Z] = R[Z]$, then $P = P^* \cap R$ is a non zero prime ideal of R. If P is of height two, then $P = rad(f, g)$ in R. But then it follows that $P^* = PK[X, Y, Z] = rad(f, g)$ in K[X, Y, Z]. If on the other hand P is of height one in R, then it is generated by an irreducible polynomial h in $R = K[X, Y]$, and if $R/(h)$ is not a Dedikend domain with torsion ideal class group, then there might be a height two prime in K[X, Y, Z] containing h and that is not a set theoretic complete intersection. This observation leads us to ask the following:

**Question 4.** R is an integral domain (not necessarily Noetherian) containing a field of characteristic zero and each prime in R[X] is a set theoretic complete intersection implies that R is of dimension one?

In the next section we answer this question in the affirmative if the coefficient ring R is Bézout.

### 3. Set theoretic complete intersection over Bézout domains

From Theorem 2, we know that (in characteristic zero case), set theoretic complete intersections of prime ideals of R[X] has more to do with the normality of R than anything else, and we now present a case further justifying this fact. To do so, we recall that an integral domain R is Bézout if each finitely generated ideal of R is principal. Our main result in this section is the following Theorem.

**Theorem 5.** *Over a finite dimensional Bézout domain R, every prime ideal of R[X] is a set theoretic complete intersection if and only if R is of dimensions one and each prime ideal of R is a set theoretic complete intersection.*

*Proof.* Follows from Proposition 6, and the proof of Theorem 7 below.

**Proposition 6.** *Let R be a finite dimensional Bézout domain. If for each maximal ideal M of R, the prime ideal MR[X] of R[X] is a set theoretic complete intersection, then R is of dimension one and each maximal ideal M of R is the radical of a principal ideal.*

*Proof.* Since R is Bézout, dimR[X] = dimR + 1, and that if P is a prime ideal of height n in R, then PR[X] is of height n in R[X]. Let M be a maximal ideal of R of height n, then MR[X] is of height n in R[X], and hence by assumption, MR[X] = $rad$ ($f_1, f_2, ..., f_n$) for some polynomials $f_1, f_2, ..., f_n$ in R[X]. Let for each i ($1 \le i \le n$), $A_{f_i}$ be the content ideal of $f_i$ in R, then as R is Bézout and each $A_{f_i}$ is finitely generated, $A_{f_i} = Ra_i$ for some $a_i$ in $A_{f_i}$. Since ($f_1, f_2, ..., f_n$) $\subseteq \sum_{i=1}^{n} A_{f_i} R[X] = \sum_{i=1}^{n} a_i R[X] \subseteq$ MR[X], it follows that the ideal ($a_1, a_2, ..., a_n$) of R is contained only in the maximal ideal M of R but not in any other maximal ideal of R. Again using the fact that R is Bézout, we have ($a_1, a_2, ..., a_n$) = $Ra$, for some $a$ in M and so M is the only maximal ideal of R containing $a$. If now P is any prime ideal of R containing $a$, then ($f_1, f_2, ..., f_n$)$\subseteq \sum_{i=1}^{n} a_i R[X] =$ $a$R[X] $\subseteq$ PR[X] $\subseteq$ MR[X] which implies that PR[X] = $rad$ ($f_1, f_2, ..., f_n$) = MR[X] in R[X] and it follows from the fact that the prime ideals of R contained in M containing the element $a$ are linearly ordered, we obtain P = M in R. Hence M = $rad(a)$ in R and of course in R[X], MR[X] = $rad(aR[X])$. Thus n = 1, and since MR[X] is a set theoretic complete intersection ideal, it is of height one. Therefore it follows that height of each maximal ideal M of R is one, and hence R is of dimension one.

**Theorem 7.** *Let R be a* Bézout domain *of Krull dimension one in which each non-zero ideal is contained in finitely many maximal ideals. Then every prime ideal of R[X] is a set theoretic complete intersection.*

*Proof.* Since R is Bézout and of dimension one, dimR[X] = 2. Let P* be any prime ideal of R[X]. Then either P* ∩ R = P is a non zero prime in R or P * ∩R = 0 in R. If P* ∩ R = P, then pick a none zero element $x$ in P and let $P_1, P_2, ... , P_n$ be the other maximal ideals of R containing $x$. Next choose an element $y$ in P not contained in any one of $P_i$ ($1 \le i \le n$), then the ideal $Rx + Ry$ of R is contained only in P but not in any other maximal ideal of R. Clearly P = $rad(Rx +Ry)$ and as R is Bézout, we have $Rx +Ry = Ra$ for some $a$ in P and so P = $rad(Ra)$. Now in R[X] we either have PR[X] = P* in which case P* = $rad$

($a$R[X]) in R[X], or else PR[X] is properly contained in P* and in that case P* is of height two, and hence by Theorem 28 of [4], P* is generated by P and some monic polynomial $f$. But then it is clear that P* = $rad$ ($a$R[X] + $f$R[X]). Next we deal with the case when P* ∩ R = 0 in R, in this case by [4, Theorem 36], we have P*K[X] = gK[X] for some polynomial g in P*, where K is the field of fractions of R. Let $A_g$ be the content ideal of g in R. Then $A_g$ = R$a$ for some $a$ in $A_g$. The remaining part of the proof is adapted from our proof of Theorem 1.2 of [3]. Let $g = b_0 + b_1X + b_2X^2 + \ldots + b_kX^k$. Then $A_g = Rb_0 + Rb_1 + \ldots + Rb_k$, and hence $a = d_0b_0 + d_1b_1 + \ldots + d_kb_k$, for some $d_0, d_1, \ldots, d_k$ in R. On the other hand for each i (0 ≤ i ≤ k), we have $b_i = c_ia$ for some $c_i$ in R. Thus $a = ac_0d_0 + ac_1d_1 + \ldots + ac_kd_k$ and from this it follows that $1 = c_0d_0 + c_1d_1 + \ldots + c_kd_k$. Therefore if $f = c_0 + c_1X + c_2X^2 + \ldots + c_kX^k$, then $g = af$ and $A_f$ = R. Hence P*K[X] = gK[X] = $a$fK[X] = fK[X] and it is clear that f is not contained in any height one prime of R[X] of the form PR[X], where P is a non-zero prime of R. Also using the fact that there is a one to one correspondence between the prime ideals of R[X] contracting to 0 in R and the prime ideals in K[X], we see that P* is the only height one prime ideal of R[X] containing f and therefore P* = $rad$(fR[X]), with this the proof is complete.

Using the second part of the above proof we have:

**Corollary 8.** *Let R be any Bézout domain. Then every prime ideal P* of R[X] contracting to zero-ideal of R is the radical of a principal ideal.*